\documentclass[preprint,a4paper,11pt]{elsarticle}

\usepackage{graphicx}
\usepackage{amsmath}
\usepackage{amscd}
\usepackage{amssymb}
\usepackage{subfig}
\usepackage{multirow}
\usepackage{float}
\usepackage{epsfig}
\usepackage{epsf}
\usepackage{stmaryrd}
\usepackage{esdiff}
\usepackage{array}
\usepackage{mathtools}
\usepackage{enumerate}
\usepackage{color}

\makeatletter
\newcommand*\bigcdot{\mathpalette\bigcdot@{.5}}
\newcommand*\bigcdot@[2]{\mathbin{\vcenter{\hbox{\scalebox{#2}{$\m@th#1\bullet$}}}}}
\makeatother

\newtheorem{thm}{Theorem}

\newtheorem{algorithm}{Algorithm}

\newtheorem{lemma}{Lemma}
\newtheorem{defn}{Definition}

\newenvironment{proof}[1][Proof]{\begin{trivlist}
\item[\hskip \labelsep {\bfseries #1}]}{\end{trivlist}}

\begin{document}

\begin{frontmatter}



\title{Network Utility Maximization Revisited:\\
Three Issues and Their Resolution}

\author[ece2]{P T Akhil\corref{cor1}}
\ead{akhilpt@gmail.com}
\cortext[cor1]{Corresponding author}

\author[ece1]{Rajesh Sundaresan}


\address[ece2]{Department of Electrical Communication Engineering, Indian Institute of Science, Bangalore 560012, India.}
 \address[ece1]{Department of Electrical Communication Engineering and Robert Bosch Centre for
Cyber Physical Systems, Indian Institute of Science, Bangalore 560012, India.}

\begin{abstract}
Distributed and iterative network utility maximization algorithms, such as the primal-dual algorithms or the network-user decomposition algorithms, often involve trajectories where the iterates may be infeasible, convergence to the optimal points of relaxed problems different from the original, or convergence to local maxima. In this paper, we highlight the three issues with iterative algorithms. We then propose a distributed and iterative algorithm that does not suffer from the three issues. In particular, we assert the feasibility of the algorithm's iterates at all times, convergence to global maximum of the given problem (rather than to global maximum of a relaxed problem), and avoidance of any associated spurious rest points of the dynamics. A benchmark algorithm due to Kelly, Maulloo and Tan (1998) [Rate control for communication networks: shadow prices, proportional fairness and stability, Journal of the Operational Research society, 49(3), 237-252] involves fast user updates coupled with slow network updates in the form of additive-increase multiplicative-decrease of suggested user flows. The proposed algorithm may be viewed as one with fast user updates and fast network updates that keeps the iterates feasible at all times. Simulations suggest that the convergence rate of the ordinary differential equation (ODE) tracked by our proposed algorithm's iterates is comparable to that of the ODE for the aforementioned benchmark algorithm.
\end{abstract}


\end{frontmatter}
\pagestyle{myheadings}
\thispagestyle{plain}
\markboth{P. T. AKHIL AND R. SUNDARESAN}{Network Utility Maximization}


\section{Introduction and the main result}
\label{sec:intro}

\subsection{Background}

We revisit the classic setting of decentralised congestion control as addressed by Kelly et al. \cite{1998xxJORS_KelMauTan}. Consider a
network with $m$ directed link resources. Let $c(l)$ be the capacity of the link $l$. There are $n$
users and each has a single fixed path. Each user sends data along its associated path with the
first vertex of the path being the source of the user's data and the last vertex being its terminus.
Let $A$ be the $m \times n$ matrix with $A_{le}=1$ if the path $e$ uses link $l$ and $A_{le}=0$
otherwise. Let $[n]$ denote the set $\{1,2,\ldots,n\}$ of users and let $[m]$ denote the set
$\{1,2,\ldots,m\}$ of links. Let $w_e:\mathbb{R}_+ \rightarrow \mathbb{R},~e \in [n]$ be the
utility functions of the users. User $e$ derives a utility $w_e(x(e))$ when sending a flow of rate
$x(e)$. The functions $w_e, e \in [n]$ are assumed to be strictly concave and increasing. Let
$x=(x(e), e \in [n])$, $w= (w_e, e \in [n])$, and $c=(c(l), l \in [m])$. Let $\mathcal{A}=\{x|x \geq
0 ~\text{and} ~ Ax \leq c\}$. Throughout, we make the standing assumption that $\mathcal{A}$ has an
interior feasible point, i.e., there exists a point for which all inequalities are
strict. The system optimal operating point solves the problem:
\begin{align}\label{eqn:maximizeSumSeparableConcaveFunctions}
 \text{System}(w,A,c):~ \max_{x \in \mathcal{A}} ~ W(x) := \sum_{e=1}^{n}w_e(x(e)).
\end{align}
The important decentralization concerns are that the network operator does not know the utility functions of the users, and the users know neither the rate choices of the other users nor the flow constraints on the network.

Kelly \cite{199701ETT_Kel} proposed the decomposition of the above problem into two subproblems, one
to be solved by each user, and the other to be solved by the network. Let $\lambda_e$ be the cost
per unit rate to user $e$ set by the network, and let $p_e$ be the price user $e$ is willing to pay.
The maximization problem solved by user $e$ is
\begin{align}\label{eqn:userob}
 \text{User}(w_e;\lambda_e):
 \max_{p_e:p_e \geq 0}~ &w_e \left(\frac{p_e}{\lambda_e}\right)-p_e .
\end{align}
If $p = (p_e, e \in [n])$ is known to the network, its optimization problem is
\begin{align}\label{eqn:maximizeSumLOgFunctions}
 \text{Network}(A,c;p):
 \max_{x\in \mathcal{A}}  ~ \sum_{e=1}^{n}p_e \log(x(e)).~~~~~~~~~~~
\end{align}
The solution to Network$(A,c;p)$ is well-known to satisfy the so-called proportional fairness
criterion: if $\mu_l, l \in [m]$ are the optimal dual price variables associated with the dual to
Network$(A,c;p)$, then
\begin{align}\label{eqn:propfairdual}
x(e)=\frac{p_e}{\sum_{l:l\in e}\mu_l}, e \in [n],
\end{align}
is the optimal solution to the network problem.
Kelly \cite{199701ETT_Kel} showed that there exist costs per unit rate $(\lambda^{\star}_e, e \in
[n])$, prices $(p^{\star}_e, e \in [n])$, and flows $(x^{\star}(e), e \in [n])$, satisfying
$p^{\star}_e = \lambda^{\star}_e \cdot x^{\star}(e)$ for $e \in [n]$ such that $p^{\star}_e$ solves
User$(w_e; \lambda^{\star}_e)$ for $e \in [n]$ and $(x^{\star}(e), e \in [n])$ solves
Network$(A,c;p^{\star})$; furthermore, $(x^{\star}(e), e \in [n])$ is the unique solution to
System$(w,A,c)$. The costs per unit rate satisfy $\lambda^{\star}_e=\sum_{l:l\in e}\mu_l^{\star}$
for some dual price variables.

In order to ensure operation at $x^{\star}$, taking the information asymmetry constraints into
account, Kelly et al. \cite{1998xxJORS_KelMauTan} proposed the following {\em fast user adaptation}
dynamics:
\begin{align}
 p_e(t)&= x(e,t) \cdot w_e^\prime(x(e,t)),~ e \in [n],\label{eqn:price}\\
 \frac{d}{dt}x(e,t)&=\kappa \cdot \left(p_e(t)-x(e,t) \cdot \sum_{l : l \in e} \mu_l(t)\right),
~e \in [n],\label{eqn:update}\\
\mu_l(t)&=\psi_l \left( \sum_{e:e \ni l} x(e,t)\right), \quad l \in [m] \label{eqn:cost},
\end{align}
where $\psi_l(y)$ is a penalty\footnote{Kelly et al. \cite{1998xxJORS_KelMauTan} suggest
two functions as examples, one of which is $\psi_l(y)=(y-c(l)+\varepsilon)^+/\varepsilon^2$ for
some $\varepsilon > 0$.} or cost per unit flow when the total flow in the link is $y$. It signifies
the level of congestion in that link. Thus $\mu_l(t)$ in (\ref{eqn:cost}) is the cost per unit flow
through link $l$, and may be interpreted as a dual variable of the network problem. The optimal
dual variables for Network$(A,c;p)$ are such that the net cost of user $e$ flow matches the price
$p_e$ paid by that user; see (\ref{eqn:propfairdual}). The network, adapts the flow
$x(e,t)$ using an additive-increase multiplicative-decrease scheme as in (\ref{eqn:update}), perhaps the first mathematical justification for the scheme already then in use in TCP/IP congestion control schemes. The network attempts to equalize, albeit slowly, the instantaneous
net cost of user $e$ flow, $x(e,t) \cdot \sum_{l : l \in e} \mu_l(t)$, to the instantaneous price
paid by that user, $p_e(t)$. On the other hand, if we differentiate (\ref{eqn:userob}) with respect
to $p_e$ and use the relation $p_e(t)=\lambda_e(t) \cdot x(e,t)$, we get that $p_e(t)$ in
(\ref{eqn:price}) maximizes User$(w_e; \lambda_e(t))$. So the users adapt {\em instantaneously} (in
comparison to the network's slower speed of adaptation) to the congestion signal. Kelly et al.
\cite{1998xxJORS_KelMauTan} provided a Lyapunov function for the dynamical system defined by
(\ref{eqn:price})-(\ref{eqn:cost}). The stable equilibrium point of the dynamical system maximizes a
relaxation of the system problem, as determined by the choice of $\psi_l(\cdot)$  in
(\ref{eqn:cost}).

The papers Kelly \cite{199701ETT_Kel} and Kelly et al. \cite{1998xxJORS_KelMauTan} are landmark papers for three reasons.
\begin{enumerate}
\item They provided perhaps the first mathematical justification for the additive-increase and multiplicative-decrease scheme then already in use for TCP/IP congestion.
\item They firmly rooted the idea of proportional fairness in the minds of network engineers.
\item They also provided the general framework to study other notions of fairness via utility functions and network utility maximization.
\end{enumerate}

\subsection{Three issues and the motivation for our work}

Despite the popularity of this approach, there are three issues we would like to highlight.
\begin{itemize}
\item $x(t)$ may not remain feasible at all times $t$.
\item $x(t)$ converges to the optimal value of a relaxation of the system problem.
\item There are multiple fixed points for the dynamics.
\end{itemize}
The first issue was highlighted in Johansson et al. \cite{2006xx_JohSolJoh}. The dynamics
(\ref{eqn:price})-(\ref{eqn:cost}) cannot then be used in systems where feasibility has to be
ensured at all times. Take for example identification of optimal flow parameters in a software defined communication network with a centralized controller. Flows may go through links of {\em fixed} capacities, and these {\em must} be respected at all times, even during the learning and exploration phases that may ensue before arrival at the final optimal flow values. This may also be required in mobile offloading settings \cite{2015xx_IosGaoHuaTas}, \cite{2018xx_NavSun} under additional assumptions of strict service level guarantees, in smart grid energy routing settings \cite{2018xxISGT_ShanEtAl}, or in road traffic settings where one simply cannot have more traffic than the road's capacity at any time, be it during the exploratory phase or otherwise.

The second issue is often circumvented via iterative algorithms where the Lagrange
multipliers or penalty functions are also adapted over time, in some examples at a slower time scale; see for example, Arrow and Hurwicz
\cite{1958xxEJES_ArrHur}, Low and Lapsley \cite{1999xx_LowLap}, Chiang et al. \cite{2007xxIEEE_ChiLowCalDoy}, Palomar and Chiang \cite{2006xxJSAC_PalChi}, and the more recent works of Gao et al. \cite{gao2013economics} and \cite{2015xx_IosGaoHuaTas}. Such approaches either assume knowledge of the utility functions at the network end or may encounter infeasible iterates, or both.

\begin{figure}[t]
\centering
\includegraphics[scale=0.35]{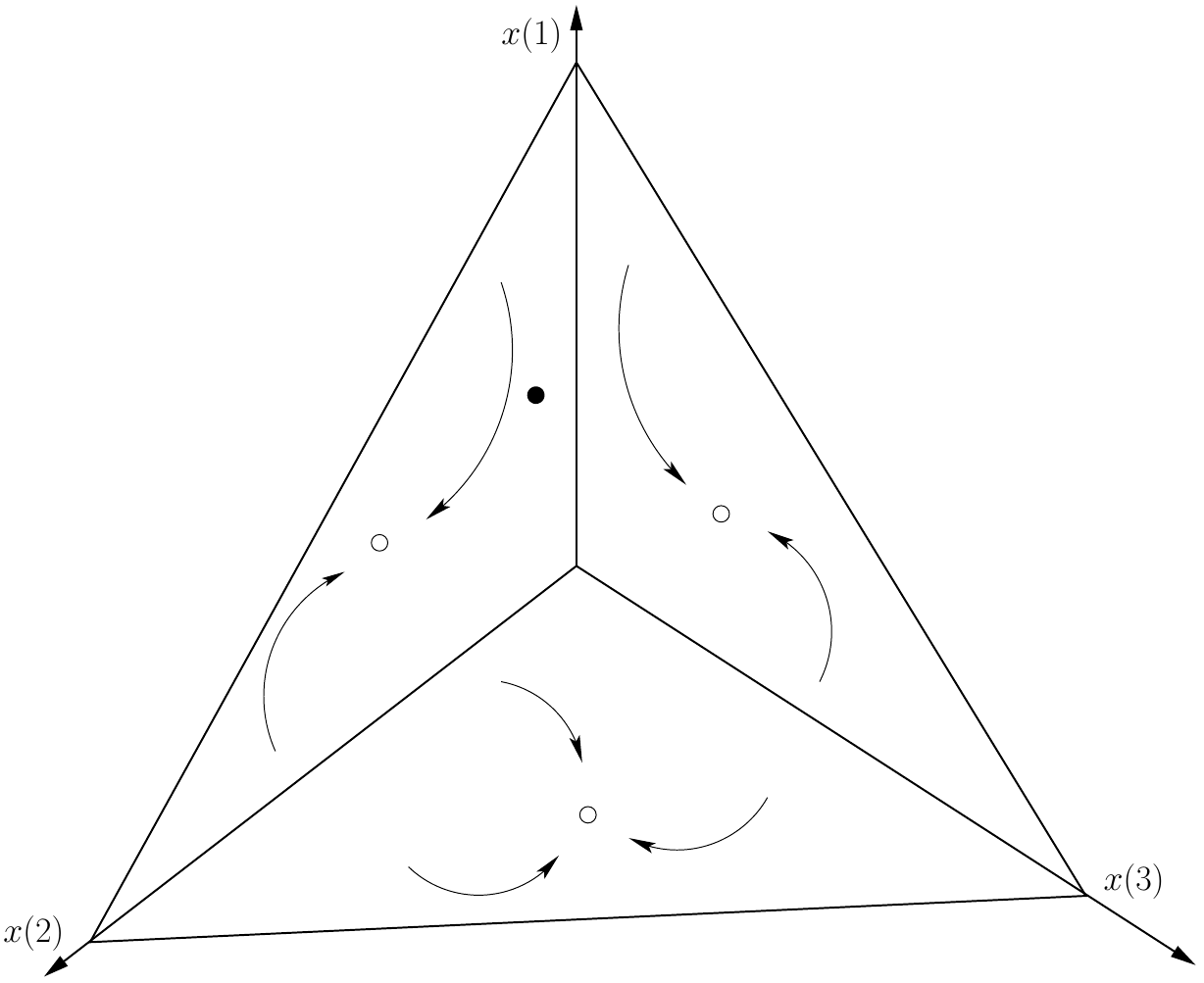}
 \caption{There are many rest points for the dynamics, with the spurious ones indicated by open circles on the facets. The solid circle is the global optimum and is on the interior. If the iteration starts on one of the facets, it remains on the facet when there is no stochasticity, and we do not have convergence to the global optimum point. Indeed, each open circle rest point is attracting for dynamics within that facet.}
 \label{fig:figure-dynamics-subspace}
 \end{figure}

The third issue is about multiple spurious rest points, other than the global optimum, for the iterative dynamics. Indeed, if $x(e,t) = 0$ for user $e$ at some time $t_0$, and if\footnote{This is the case when the marginal utility for user $e$ at supplied rate $r=0$, $w'_e(0+)$, is finite. The quantity $\lim_{r \rightarrow 0} r w'_e(r)$ can be nonzero only when $w'_e(0+) = \infty$, for e.g., $w_e(r) = \log r$.} $\lim_{r \rightarrow 0} rw'_e(r) = 0$, then from (\ref{eqn:price}), the user's willingness to pay $p_e(t_0) = 0$, and from (\ref{eqn:update}) one gets $x(e,t) \equiv 0, t \geq t_0$, i.e., the iterates never exit the facet defined by $x_e = 0$. See Figure \ref{fig:figure-dynamics-subspace}. When there is no stochasticity, there is no exit from this facet, and the iterates converge to a rest point for the dynamics different from the global optimum. Further, there is no a priori guarantee that these rest points are not attracting, and so there may be no exit even under stochasticity if the iterates start sufficiently close to these rest points.

The purpose of this paper is to provide an algorithm that circumvents these three issues.

\subsection{Most relevant related works}

The literature on network utility maximization is so vast that we will not be able to do justice to twenty years of literature on the topic. However, we will focus on works where the iterates remain feasible at all times.  There are three works, Hochbaum \cite{1994xxMOR_Hoc}, Mo and Walrand  \cite{2000xxTON_MoWal} and La and Anantharam \cite{2002xxTON_LaAna}, that are very relevant to our contribution which we bring to the reader's attention. A greedy algorithm proposed by Hochbaum \cite{1994xxMOR_Hoc} can be adapted to solve the system problem with iterates remaining feasible at all times and without full knowledge of the utility functions at the network side. Though the algorithm circumvents the issues highlighted above, it works only when the  set of feasible flows forms the ``independent set of a
polymatroid". This is the case when the network has, for example, a single source and multiple sinks
or when the network has multiple sources but a single sink.

Mo and Walrand \cite{2000xxTON_MoWal} proposed a window-based rate control mechanism that converges
to the solution to Network$(A,c;p)$ for a fixed $p$. The window update rule of
\cite{2000xxTON_MoWal} uses only delay information provided to the user (propagation and
round-trip delays).

La and Anantharam \cite{2002xxTON_LaAna}
proposed two algorithms that solve the system problem using the decomposition of Kelly et al.
\cite{1998xxJORS_KelMauTan}. The first algorithm incorporates the solution to the user problem into
the window update rule of \cite{2000xxTON_MoWal}. The second algorithm of La and Anantharam
\cite{2002xxTON_LaAna} explicitly finds the solution to the user problem and the network problem in
each iteration. Although their simulations showed the convergence of the algorithm for general
networks, a rigorous proof was given only for the case of a network with a single link. Their
algorithm additionally imposes more stringent conditions on the utility functions than those assumed
in this paper.

\subsection{Our main result} In this paper, we propose a discrete-time algorithm (see Algorithm \ref{alg:mainalgorithm} below)
that (1) remains feasible at all times, (2) converges to the desired global maximum of System$(w,A,c)$ (rather than to global maximum of a relaxed problem), and (3) therefore avoids spurious traps. The corresponding continuous time dynamics also shares the same properties. In comparison to \cite{2002xxTON_LaAna}, our algorithm applies to more general networks and a larger class of objective functions.

We now set up the notation to describe the algorithm. For a set of flows, abusing notation, write $p_e(x(e)):=x(e) \cdot w_e^{\prime}(x(e))$ as per
(\ref{eqn:price}), and set $p(x) = (p_e(x(e)), e \in [n])$. Write
\begin{align}
  \label{eqn:functionT}
  T(x) \coloneqq \arg \max_{y \in \mathcal{A}} \sum_{e=1}^{n}p_e(x(e)) \log(y(e))
\end{align}
for the solution to Network$(A,c;p(x))$. If $p_e(x(e))=0$ for some $e$, then the objective function
in (\ref{eqn:functionT}) is not strictly concave over $\mathcal{A}$. The optimization problem
(\ref{eqn:functionT}) may then have multiple solutions, and so $T(x)$ is to be viewed as a
set-valued mapping whose values are convex and compact subsets of $\mathcal{A}$. Define
$a_k\coloneqq \frac{1}{k+1}, k =0,1,2,\ldots$.

\begin{algorithm}\label{alg:mainalgorithm}
\begin{enumerate}
\item Initialize $x^{(0)} \in \mathcal{A}$  such that $x^{(0)}(e) > 0,~e \in [n]$. Initialize $k=0$.
\item User update:\label{step:step2}
 \begin{align}\label{eqn:userupdate}
  p_e^{(k)}=x^{(k)}(e) \cdot w_e^{\prime}(x^{(k)}(e)),~e \in [n].
 \end{align}
\item Network update:\label{step:step4}

Find a point $v \in T(x^{(k)})$ and set:
\begin{align}x^{(k+1)}= x^{(k)}+a_{k+1}(v-x^{(k)}).\label{eqn:stocapproxscheme}
\end{align}
\item Set $ k\leftarrow k+1$ and go to step \ref{step:step2}.
\end{enumerate}
\end{algorithm}

Our main result is the following theorem.
\vspace{2mm}
\begin{thm}
\label{thm:maintheorem}
Assume that $\mathcal{A}$ has an interior feasible point. The iterates $x^{(k)}$ of Algorithm \ref{alg:mainalgorithm} converge to $x^{\star}$, the optimal solution to the system problem, i.e.,
$x^{(k)} \rightarrow x^{\star}$ as $k \rightarrow \infty$.
\end{thm}

\subsection{The three issues are now resolved}

Recall that our objective is to address the three main issues in the dynamics of (\ref{eqn:price})-(\ref{eqn:cost}), viz., the non-feasibility of the iterates, their convergence to a different solution -- the solution to some relaxed problem, or their convergence to local maxima traps on one of the facets. We now argue that these issues disappear for Algorithm \ref{alg:mainalgorithm}.

Observe that, in Algorithm \ref{alg:mainalgorithm}, the users exhibit the same fast adaptation as in the dynamics (\ref{eqn:price}). But in the network update, iterate $x^{(k+1)}$ is a convex combination of $x^{(k)}$ and $v \in T(x^{(k)})$ which, by induction, remains in the feasible set for all $k$. This resolves the feasibility issue that plagues the dynamics (\ref{eqn:price})-(\ref{eqn:cost}).
In the proof we will argue that the iterates track the differential inclusion
\begin{align}\label{eqn:mainalgdi}
\frac{d}{dt}x(t) \in T(x(t)) - x(t);
\end{align}
we will in fact see that the solution to this differential inclusion also remains feasible at all times.

Theorem \ref{thm:maintheorem} asserts that the iterates converge to the global optimum of the system problem. This resolves the issue that the dynamics (\ref{eqn:price})-(\ref{eqn:cost}) converge to the solution to a relaxed problem different from the original system problem.

\begin{figure}[t]
\centering
\includegraphics[scale=0.35]{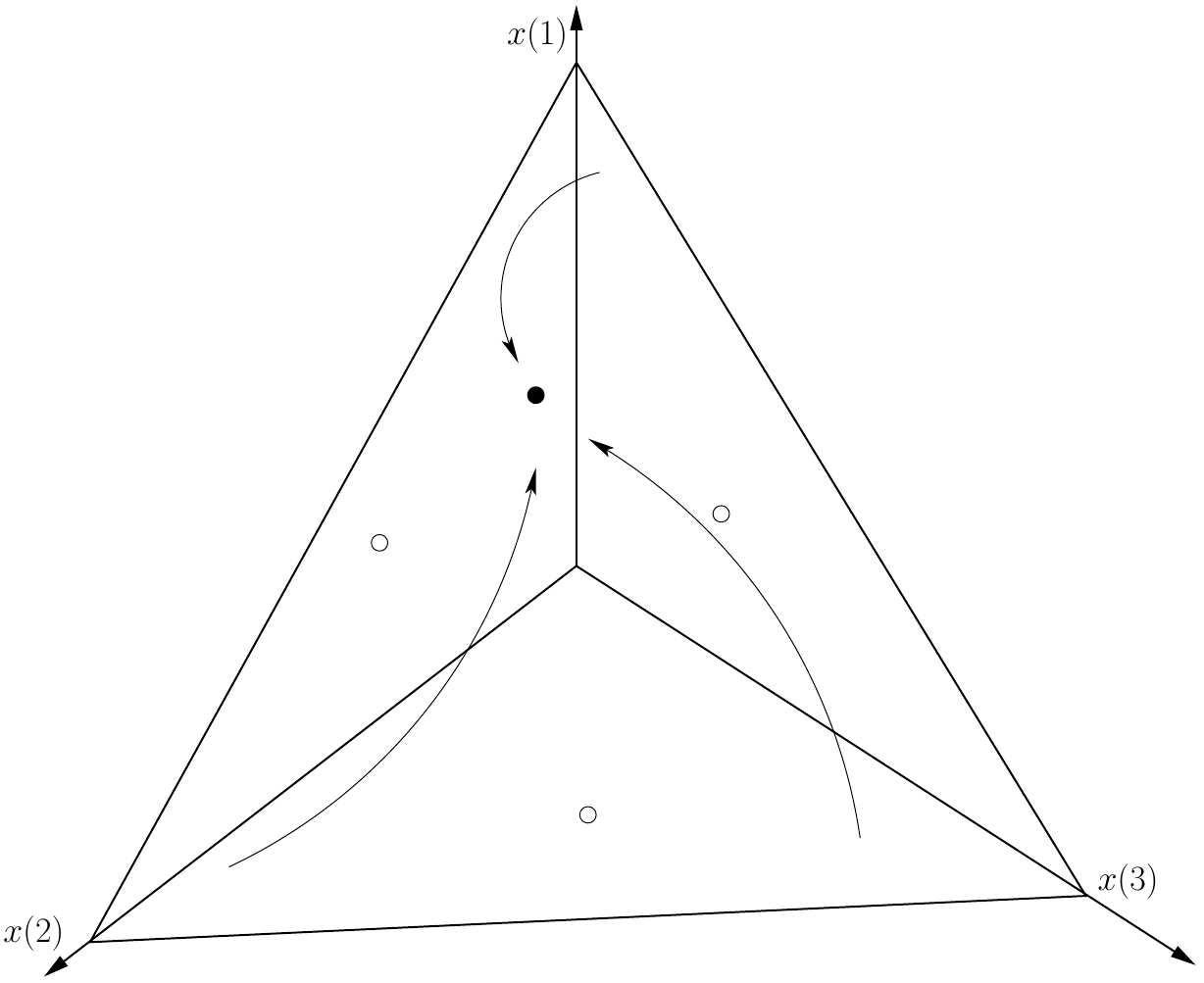}
 \caption{Dynamics of Algorithm 1. Even if we start close to the facet, so long as the initial point is on the interior, there is convergence to the global optimum.}
 \label{fig:figure-dynamics-optimum}
 \end{figure}

The assertion that there is convergence to the global optimum resolves the third issue as well of avoidance of spurious rest points on the facets. This is particularly interesting since there is no stochasticity in our algorithm. See Figure \ref{fig:figure-dynamics-optimum}. We must however start in the interior, but any arbitrary interior feasible point will work.

\subsection{Nontriviality of our contribution}
1) {\em Slowdown is essential}: It is worthwhile to ask whether the convergence of Algorithm \ref{alg:mainalgorithm} could be sped up using a constant step size $a_k \equiv \gamma,$ where $0 <\gamma < 1,$ in the network update. In fact, Hou and Kumar \cite{2010xxINF_HonKum} proposed a variant of Algorithm \ref{alg:mainalgorithm} with constant step size in the network update rule. This was done in the context of delay constrained throughput maximization in wireless networks. Step sizes determine the sizes of exploration. Larger step sizes involve faster exploration but also increased variance. In  \ref{sec:stepsizecounterexample}, we provide a counterexample to show that the algorithm of Hou and Kumar with a constant step size \emph{may not converge}. The choice of step sizes is therefore a delicate matter. A sufficient condition on step sizes $a_k, k\geq 1$, for the algorithm to converge is the usual conditions typical in stochastic approximation literature.
\begin{align}
\lim_{k \rightarrow \infty}a_k = 0,\mbox{ and} ~\sum_{k=1}^{\infty}a_k \rightarrow \infty.
\end{align}
The second condition ensures that there is enough exploration while the first condition gradually brings the exploration or learning rate to $0$.

2) {\em Technical issues in the proof of convergence}: The main technical issues to surmount in showing the convergence of Algorithm \ref{alg:mainalgorithm} are (a) the dynamics in (\ref{eqn:mainalgdi}) also have
multiple fixed points and it is \emph{nontrivial} to show convergence to the global optimum; (b) $T(x)$ is not necessarily a continuous function of $x$; see Section \ref{subsec:NeedForDI} and \ref{appsec:DiscontinuousT}. One must then study a differential inclusion, and the solution space may not be unique in general. An additional issue is the lack of stochasticity. We must then show that local maxima traps can yet be escaped. The main technical contribution is that these issues can be surmounted at least in this NUM problem with sufficient convexity structure.

3) {\em Wide applicability}: Iterative procedures for network utility maximization have wide applicability. Road traffic control, software defined network controllers with strict service level guarantees, offloading of data-traffic to WiFi providers and femto cell providers (\cite{2015xx_IosGaoHuaTas} and \cite{2018xx_NavSun}) with strict call handling requirements, and smart grid energy routing \cite{2018xxISGT_ShanEtAl} are all settings where feasibility should be met at all times and where convergence to global optimality is desirable. Our algorithm is applicable in all these settings.

4) {\em Decentralized implementation}: At first glance, Algorithm \ref{alg:mainalgorithm} {\em appears} to need a central entity that computes the solution to Network$(A,c;p)$, i.e. $T(x)$, and Section \ref{sec:complexity} describes algorithms to compute $T(x)$ efficiently for some class of networks. However, the central entity is {\em not essential} because we can use Mo and Walrand's algorithm \cite{2000xxTON_MoWal} to find $T(x)$ for a fixed $p$; that algorithm uses only the information available at the user end. The $p$ can then be adapted (user updates) at a slower time scale. This enables the potential use of Algorithm \ref{alg:mainalgorithm} in a {\em distributed setting} in large scale networks.

\subsection{Organization}
The rest of the paper is organized as follows. In Section \ref{sec:convergence}, we prove Theorem
\ref{thm:maintheorem}. In Section \ref{sec:complexity}, we address the complexity of identifying the
proportionally fair solution point for the network problem. We provide an example of a network where
flows aggregate into a `main branch', reminiscent of traffic from the suburbs flowing into an
arterial highway leading to the downtown of a large city, for which the complexity to solve the
network problem is $\mathcal{O}(n)$. We also argue that this complexity is manageable
($\mathcal{O}(n \log n)$ plus computations for feasibility checks) in situations where the feasible
set is a polymatroid, for example, when all flows either originate or terminate at a single vertex.
We also see in simulations that the dynamics in (\ref{eqn:mainalgdi}) converge to the
equilibrium  at a faster rate than the dynamics of (\ref{eqn:price})-(\ref{eqn:cost}) for identical
speed parameters $\kappa$. In Section \ref{sec:conclusion}, we end the paper with some
concluding remarks.

\section{Proof of Convergence}
\label{sec:convergence}
The update equation in step \ref{step:step4} of Algorithm \ref{alg:mainalgorithm} is a standard
stochastic approximation scheme but without the stochasticity. A common method to analyze the
asymptotic behavior of such schemes is the dynamical systems approach based on the theory of
ordinary differential equations (ODE). But $T(x)$ being a set valued map necessitates the use of
differential inclusions.

The outline of the proof is as follows. We will first characterize the fixed points of the
mapping $T$. We will then argue that the system optimal point is one of the finitely many fixed
points of the mapping $T$. We will next show that the solution to the differential inclusion in
(\ref{eqn:mainalgdi}) models the asymptotic behavior of the iterates $x^{(k)}$.
Following this, We will show that every solution to the differential inclusion converges to
one of the fixed points of $T$ via Lyapunov theory. Finally, though there may be many fixed points,
we will prove that the fixed point to which the solution to the differential inclusion converges
as $t \rightarrow \infty$ is the system optimal point.

\subsection{Characterization of the fixed points of $T(x)$}
\begin{defn}
A point $x$ is a fixed point of the set valued map $T$ if $x \in T(x)$.
\end{defn}

Let $S \subset [n]$. Let $\mathcal{A}|_{S}$ be the subset of $\mathcal{A}$ whose points have support contained within $S$. Define a subproblem of the system problem as
\begin{align}
\text{Subsystem}(w,A,c,S):
\max_{y \in \mathcal{A}|_{S}}~&~\sum_{e\in
S}w_e(y(e)).~~~~~~~\label{eqn:objsubproblem}~~~~~~~~~~~~~~~~~~
\end{align}

\begin{lemma}\label{lem:finitefixedpoint}
Let $\overline{x}$ be a fixed point of the mapping $T$.  Let $S=\{e:\overline{x}(e) > 0\}$. Then
$\overline{x}$ is the unique optimal solution to the Subsystem$(w,A,c,S)$.
\end{lemma}
\begin{proof}
We have $\overline{x}(e)=0$ for all $e \in S^{c}$. If $\lim_{x(e) \downarrow
0}x(e)\cdot w_{e}^{\prime}(x(e))>0$ for some $e \in S^{c}$, then any element  $y \in T(\overline{x})$ has $y(e)>0$ which contradicts the fact
that $\overline{x}$ is a fixed point. Hence $p_e(\overline{x}(e))=\overline{x}(e)\cdot w_{e}^{\prime}(\overline{x}(e))=0$ for all $e \in S^{c}$,
and we may write
\begin{align}
\overline{x} \in T(\overline{x}) = \arg\max_{y \in \mathcal{A}} ~~\sum_{e\in S}p_e(\overline{x}(e)) \log(y(e)).\label{eqn:subnetprob}
\end{align}
Since $\overline{x} \in \mathcal{A}|_S \subset \mathcal{A}$, we also have that  $\overline{x}$ maximizes (\ref{eqn:subnetprob}) over
$\mathcal{A}|_{S}$, i.e.,
\begin{align}
\overline{x} = \arg\max_{y \in \mathcal{A}|_S}~~\sum_{e\in S}p_e(\overline{x}(e))
\log(y(e)).\label{eqn:subnetprob2}
\end{align}
Let $\mu_l,l \in [m]$ and $\eta_e, e \in [n]$ be the optimal dual variables for the network
subproblem (\ref{eqn:subnetprob2}). Its Karush-Kuhn-Tucker (KKT) conditions are
\begin{align}
&\frac{p_e(\overline{x}(e))}{\overline{x}(e)}=\sum_{l:l\in e}\mu_l - \eta_e,~e \in S, \label{eqn:lagderivsub}\\
 &\mu_l\cdot \left (\sum_{e:e \ni l}\overline{x}(e)-c(l)\right)=0,~l\in [m],\label{eqn:compslackmusub}\\
& \eta_e \cdot \overline{x}(e)=0,~e \in [n], \label{eqn:compslacketasub}\\
 &\eta_e \geq 0, ~e \in S, ~\mu_l \geq 0,~l\in [m] ~\text{and}~  \overline{x}\in \mathcal{A}|_{S}.\label{eqn:feaspriduasub}
\end{align}
Since $\frac{p_e(\overline{x}(e))}{\overline{x}(e)}=w_e^\prime(\overline{x}(e))$, it is easy to see
that equations (\ref{eqn:lagderivsub}-\ref{eqn:feaspriduasub}) are the KKT conditions of
Subsystem$(w,A,c,S)$ as
well. Since $\overline{x}/2$ is an interior feasible point of $\mathcal{A}|_S$, KKT conditions are
sufficient
for optimality in (\ref{eqn:objsubproblem}) and $\overline{x}$ is the optimal solution to
Subsystem$(w,A,c,S)$. Uniqueness follows from the strict concavity of $w_e, e \in S$. $\hfill \blacksquare$
\end{proof}
Observe that there are only a finite number of sub-problems of the form Subsystem$(w,A,c,S)$, $S \subset [n]$. As a consequence of Lemma \ref{lem:finitefixedpoint}, every fixed point of $T$ is the unique optimal solution to Subsystem$(w,A,c,S)$ for some set $S \subset [n]$. Hence there are only finitely many fixed points of $T$, each corresponding to a sub-problem Subsystem$(w,A,c,S)$ for some $S \subset [n]$.

Is every solution to Subsystem$(w,A,c,S)$ a fixed point of the mapping $T$? The possibility that $\lim_{x(e) \downarrow 0}x(e)\cdot w_{e}^{\prime}(x(e))>0$ for an $e \in S^{c}$ and the first step of the proof of Lemma \ref{lem:finitefixedpoint} says this is not always true. However, we can assert the following.

\begin{lemma}
The global maximum of the system problem, $x^{\star}$, is a fixed point of the mapping $T$.
\end{lemma}
\begin{proof}
$x^{\star}$ solves the system problem. Let $S=\{e:x^{\star}(e) >0\}$. $S$
can be a proper subset of $[n]$. Let $\mu_l,l \in [m]$ and $\eta_e, e \in [n]$ be the optimal dual
variables of the system problem. We then have
\begin{align}
&w_e^{\prime}(x^{\star}(e))=\sum_{l:l\in e}\mu_l - \eta_e,~e \in S, \label{eqn:lagderivsys}\\
&w_e^{\prime}(0)=\sum_{l:l\in e}\mu_l - \eta_e,~e \in S^{c}, \label{eqn:lagderivsys2}\\
&\mu_l\cdot \left(\sum_{e:e\ni l}x^{\star}(e)-c(l)\right)=0,~l\in
[m],\label{eqn:compslackmusys}\\
 &\eta_e \cdot x^{\star}(e)=0, ~e \in [n], \label{eqn:compslacketasys}\\
&\eta_e \geq 0,~ e \in [n],~ \mu_l \geq 0,~l\in [m] ~\text{and}~ x^{\star} \in \mathcal{A}.
\label{eqn:feaspriduasys}
\end{align}
 Observe that $w_e^{\prime}(0)$ is finite for an $e \in S^{c}$; otherwise a small increase in
$x^{\star}(e)$ and a corresponding decrease in $x^{\star}(i)$ for a suitable $i \in S$ (which has
finite $w_i^{\prime}(x^{\star}(i))$)  will result in a feasible flow that has a larger objective
function value. Hence $p_e(x^{\star}(e))=0$ for $e \in S^{c}$. Since
$w_e^{\prime}(x^{\star}(e))=\frac{p_e(x^{\star}(e))}{x^{\star}(e)}$ for all $e \in S$, it follows
from (\ref{eqn:lagderivsys})-(\ref{eqn:feaspriduasys}) that $x^{\star},(\mu_l, l\in [m]),
(\tilde{\eta}_e=\eta_e, e \in S)$ and $(\tilde{\eta_e}=\eta_e+w_e^{\prime}(0), e \in S^{c})$
satisfy the KKT conditions of the problem (\ref{eqn:functionT}). Hence $x^{\star} \in T(x^{\star})$. $\hfill \blacksquare$
\end{proof}

The above result provides a motivation to search for the global maximum by setting up a dynamics that will converge to a fixed point of $T$.

\subsection{Need for the theory of differential inclusions}\label{subsec:NeedForDI}
We now describe the issues that make it necessary to use differential inclusions to study the
asymptotic behavior of $x^{(k)}$. $T(x)$ is the set of points that solve (\ref{eqn:functionT}). If
$p_e(x(e)) = 0$ for some $e$ at a point $x \in \mathcal{A}$, then the objective function in
(\ref{eqn:functionT}) is not strictly concave. Hence there can be multiple points that solve
(\ref{eqn:functionT}). A continuous selection\footnote{A continuous selection $f$ from the set-valued map $T$ is a continuous function $f:\mathcal{A}\rightarrow \mathcal{A}$  with $f(x) \in T(x)$ for each $x \in \mathcal{A}$.} from $T(x)$ allows the use of differential equations
to analyze the stochastic approximation scheme in (\ref{eqn:stocapproxscheme}). A natural question
that arises is whether there is such a continuous selection from $T(x)$. We give an example in \ref{appsec:DiscontinuousT} showing that such a selection is not always possible.
\subsection{Differential Inclusions: Preliminaries}
In this section, we define a differential inclusion and state relevant results from \cite{2005xxSIAM_BenHofSor} that are used to show the
convergence of Algorithm \ref{alg:mainalgorithm}.
Let $F: \mathbb{R}^{n} \rightarrow \mathbb{R}^{n}$ be a set valued map. Consider the following differential inclusion:
\begin{align}\label{eqn:DI}
\frac{dx}{dt}\in F(x).
\end{align}
A solution to the differential inclusion in (\ref{eqn:DI}) with initial condition $x_{0} \in
\mathbb{R}^{n}$ is an absolutely continuous function $x:\mathbb{R} \rightarrow \mathbb{R}^{n}$ that
satisfies (\ref{eqn:DI}) for almost every $t \in \mathbb{R}$. The following conditions are
sufficient for the existence of a solution to the differential inclusion (\ref{eqn:DI}):
\begin{enumerate}
\item \label{itm:property1}$F(x)$ is nonempty, convex and compact for each $x \in \mathbb{R}^{n}$.
\item \label{itm:property2}$F$ has a closed graph.
\item \label{itm:property3}For some $K>0$, for all $x \in \mathbb{R}^{n}$, $F$ satisfies the following condition
\begin{align}\label{eqn:Fboundedness}
\sup_{z \in F(x)} \vert\vert z\vert\vert \leq K (1+\vert\vert x\vert \vert).
\end{align}
\end{enumerate}
The stochastic approximation scheme with iterates in $\mathbb{R}^{n}$ is given as
\begin{align}
y^{k+1} \in y^{(k)}+a_{k+1}(F(y^{(k)}) +U^{(k+1)}),\label{eqn:stocapproxschemenoise}
\end{align}
where $a_k$ satisfy the usual conditions:
\begin{align}
\lim_{k \rightarrow \infty}a_k = 0, ~\sum_{k=1}^{\infty}a_k \rightarrow \infty,
\end{align}
and $U^{(k)} \in \mathbb{R}^{n}$ are deterministic or random perturbations.

Let $t(0)=0, t(k)=\sum_{i=1}^{k}a_i$. Let $r_{y}:\mathbb{R}_{+}\rightarrow \mathbb{R}^{n}$ be a
continuous piece-wise linear function formed by the interpolation of $y^{(k)}$ as in
\begin{equation}
\label{eqn:linintpol}
r_{y}(t) = y^{(k)}+\frac{y^{(k+1)}-y^{(k)}}{t(k+1)-t(k)}\cdot (t-t(k)),\quad \forall ~t \in [t(k),t(k+1)).
\end{equation}
\begin{defn}
(A perturbed solution to (\ref{eqn:DI})). Let $U:\mathbb{R}_{+} \rightarrow \mathbb{R}^n$ be locally integrable function such that
\begin{align*}
\lim_{t\rightarrow \infty} \sup_{0 \leq v \leq T}\left\vert\left\vert\int_{t}^{v}U(s)ds\right\vert\right\vert = 0.
\end{align*}
Let $\delta:[0,\infty) \rightarrow \mathbb{R}$ be a function such that $\delta(t) \rightarrow 0$ as $t \rightarrow \infty$. Define
\begin{align}
F^{\epsilon}(x) = \{y \in \mathbb{R}^{m}: \exists z: \vert\vert z-x\vert\vert < \epsilon, d(y,F(z)) < \epsilon\},
\end{align}
where $d(y,F(z)) = \inf\{\vert\vert y-q\vert\vert:q \in F(z)\}$.
An absolutely continuous function $y:[0, \infty)\rightarrow \mathbb{R}^{n}$ is a perturbed solution to (\ref{eqn:DI}) if there exists
$U:\mathbb{R}_+ \rightarrow \mathbb{R}^{n}$ and $\delta:\mathbb{R}_+ \rightarrow \mathbb{R}^{n}$ as above such that
\begin{align}
\frac{dy}{dt} - U(t) \in F^{\delta(t)}(y(t)).
\end{align}
for almost every $t \in \mathbb{R}_{+}$.
\end{defn}
The following lemma, taken from \cite{2005xxSIAM_BenHofSor}, gives conditions on $y^{(k)}$ and $U^{(k)}$ for $r_y$ to be a perturbed solution to
(\ref{eqn:DI}).
\begin{lemma}\label{lem:conditionforps}
\cite[Prop. 1.3]{2005xxSIAM_BenHofSor} Suppose $y^{(k)}$ is bounded, i.e., $\sup_{k} \vert \vert y^{(k)} \vert \vert < M < \infty$, and for all
$T > 0$,
\begin{align}
\lim_{s \rightarrow \infty} \sup \Biggl\{ \left\vert\left\vert\sum_{k=s}^{i-1} a_{k+1}
U^{(k+1)}\right\vert\right\vert:i=s+1,s+2, \ldots, m(t(s)+T)
\Biggr \}=0, \label{eqn:conditionforps2}
\end{align}
where $m(t) = \sup\{m:t(m) \leq t\}$. Then $r_{y}(t)$ is a perturbed solution of the differential
inclusion (\ref{eqn:DI}).
\end{lemma}
\begin{defn}
A compact set $L$ is an internally chain transitive set if for any $x, y \in L$ and every $\epsilon
> 0 , T > 0$, there exists $l \in \mathbb{N}$, solutions $x_1,x_2,\ldots,x_l$ to (\ref{eqn:DI}) and
$t_i > T,~ \forall ~ i$, that satisfy the following.
\begin{enumerate}
\item $x_i(t) \in L$ for all $ 0 \leq t \leq t_i$ and for all $i \in [l]$,
\item $\vert\vert x_i(t_i)- x_{i+1}(0) \vert \vert \leq \epsilon $ for all $i \in [l-1]$,
\item $\vert\vert x_1(0)- x \vert \vert \leq \epsilon$ and $\vert\vert x_l(t_l)- y\vert \vert \leq \epsilon$.
\end{enumerate}
We shall call the sequence $(x_1,x_2,\ldots,x_l)$ as an $(\epsilon,T)$ chain in $L$ from $x$ to $y$.
\end{defn}
The following lemma, again taken from \cite{2005xxSIAM_BenHofSor}, characterizes the limit set of
a perturbed solution.
\begin{lemma}\cite[Thm. 3.6]{2005xxSIAM_BenHofSor}\label{lem:psisict}
Let $r$ be a perturbed solution to (\ref{eqn:DI}). Then the limit set of $r(\cdot)$ $L(r)\coloneqq \bigcap_{t \geq 0}\{ r(s) : s \geq t\}$ is internally chain transitive.
\end{lemma}
\subsection{Convergence analysis}
We proceed to prove the convergence of $x^{(k)}$, the iterates put out by Algorithm \ref{alg:mainalgorithm}, to the optimal solution to the system problem. Observe that $T$ maps points in $\mathcal{A}$ to
itself. We show that $x^{(k)}$ asymptotically tracks the solution to the differential inclusion
\begin{align}\label{eqn:DI2}
\frac{dx}{dt}\in F(x),
\end{align}
where
\begin{align} \label{eqn:Fdefinition}
F(x) \coloneqq T(P_{\mathcal{A}}(x))-x,
\end{align}
for $x \in \mathbb{R}^{n}$. $P_{\mathcal{A}}(x)$ is the projection of $x$ onto the set $\mathcal{A}$. We then find a Lyapunov function for the dynamics in (\ref{eqn:DI2}) to show its convergence to the system optimal point $x^{\star}$.
The following lemma establishes that the set-valued map $F$ defined in (\ref{eqn:Fdefinition}) has some good properties; it turns out that these are sufficient for the differential inclusion (\ref{eqn:DI}) to have a solution.

\begin{lemma}\label{lem:DIsolnexists}
For each $x \in \mathbb{R}^{n}$, $F(x)$ is nonempty, convex and compact. Furthermore, $F$ has the closed graph property and satisfies
(\ref{eqn:Fboundedness}).
\end{lemma}
\begin{proof}
The objective function of the network problem is continuous and the constraint set $\mathcal{A}$ is
compact. The maximum exists due to the Weierstrass theorem. Also, the set of maximizers is closed
and convex. Thus $T(P_{\mathcal{A}}(x))$ is nonempty, convex and compact, and hence so is $F(x)$.

We next prove the closed graph property of $F$. A function has the closed graph property if it is
upper hemicontinuous. The objective function in (\ref{eqn:functionT}), $\sum_{e=1}^{n}p_e(x(e))\log
(y(e))$, is jointly continuous\footnote{Here we may take $\log(y(e))$ to be continuous at $y(e) =
0$ with $\log \, 0 \coloneqq -\infty$ because then $\lim_{x(e) \downarrow 0} \log(y(e))=\log\, 0$
and this sequential continuity is all that is needed to apply Berge's maximum theorem
\cite[p. 116]{1963xxTS_Ber}} in $x$ and $y$. Also, the constraint set of the network problem does
not vary with $x$. By Berge's maximum theorem \cite[p. 116]{1963xxTS_Ber}, $T(x)$ is upper
hemicontinuous. Since $P_{\mathcal{A}}(x)$, the projection onto the  convex set $\mathcal{A}$, is
continuous, the composition $T(P_{\mathcal{A}}(x))$ is upper hemicontinuous. Consequently, $F(x)$ is
upper hemicontinuous and hence has the closed graph property.
Finally,
\begin{align}
\sup_{z \in F(x)} \vert\vert z\vert\vert &= \sup_{z \in T(P_{\mathcal{A}}(x))} \vert\vert z-x \vert\vert \leq \sup_{z \in \mathcal{A}}\vert\vert z-x \vert\vert \leq \sup_{z \in \mathcal{A}} \vert\vert z \vert\vert + \vert\vert x
\vert\vert \nonumber \\& \leq K + \vert\vert x \vert\vert \leq K (1 + \vert\vert x \vert\vert),\text{~where~} K > \max\{\sup_{z \in \mathcal{A}}\vert\vert z \vert \vert,1\}.\nonumber
\end{align}
\end{proof}
\begin{lemma}\label{lem:rxisps}
Let $r_x(t)$ be obtained by the linear interpolation of $x^{(k)}$ as given in (\ref{eqn:linintpol}).
Then $r_x(t)$ is a perturbed solution to the differential inclusion (\ref{eqn:DI}) with $F$ defined
as in (\ref{eqn:Fdefinition}).
\end{lemma}
\begin{proof}
We first show that $x^{(k)} \in \mathcal{A}$. Observe that $x^{(0)} \in \mathcal{A}$. Assume $x^{(k-1)} \in \mathcal{A}$. Since $T(x^{(k-1)}) \in \mathcal{A}$ and $x^{(k)}$ is a convex combination of $x^{(k-1)}$ and $T(x^{(k-1)})$, we have $x^{(k)} \in \mathcal{A}$. It follows that
\begin{align*}
F(x^{k}) = T(P_{\mathcal{A}}(x^{(k)}))-x^{(k)} = T(x^{(k)})-x^{(k)} .
\end{align*}
We now see that the update equation in (\ref{eqn:stocapproxscheme}) is the same as the stochastic
approximation scheme in (\ref{eqn:stocapproxschemenoise}) with $U^{(k)}= 0$ for all $k$.
Observe that $x^{(k)}$ is bounded because $x^{(k)} \in \mathcal{A}$ for all $k$ and $\mathcal{A}$
is compact; since $U^{(k+1)}= 0$, the condition in (\ref{eqn:conditionforps2}) is trivially
satisfied.
Hence, by Lemma \ref{lem:conditionforps}, $r_x(t)$ is a  perturbed solution. $\hfill \blacksquare$
\end{proof}
\vspace{1mm}
We restrict our attention to solutions of (\ref{eqn:DI}) with initial condition $x(0) \in \mathcal{A}$. Since $T(x) \in \mathcal{A}$, $x(t)$ lies
in $\mathcal{A}$ for all $t$.
Define \[\Phi_{t}(x_0)\coloneqq \{x(t):x ~\text{solves}~ (\ref{eqn:DI}),~x(0)=x_0\}.\]
\begin{defn}
Let $\Lambda$ be a subset of $\mathcal{A}$. Let $V:\mathcal{A} \rightarrow \mathbb{R}$ be a continuous function such that $V(y) < V(x),~ y \in \Phi_t(x), x \in \mathcal{A}\backslash \Lambda$ and $V(y) \leq V(x),~ y \in \Phi_t(x), x \in \Lambda.$ Then $V$ is called a Lyapunov function for $\Lambda$.
\end{defn}
Define $V:\mathcal{A} \rightarrow \mathbb{R}$ as
\begin{align}\label{eqn:objlyapunov}
V(x)\coloneqq \sum_{e=1}^{n} w_{e}(x^{\star}(e))-\sum_{e=1}^{n}w_{e}(x(e)).
\end{align}
\begin{lemma}\label{lem:Vlyapunov}
Let $\Lambda$ be the set of fixed points of $T$. The function $V$ in (\ref{eqn:objlyapunov}) is a
Lyapunov function for $\Lambda$.
\end{lemma}
\begin{proof}
Let $x \in \mathcal{A}$ and $v \in T(x)$. We have, from the definition of $T(x)$ in
(\ref{eqn:functionT}), that
\begin{align}\label{eqn:strictinequality}
\sum_{e=1}^{n}p_e(x(e)) \log (v(e)) \geq \sum_{e=1}^{n} p_e(x(e)) \log (x(e))
\end{align}
because $v \in T(x)$ maximizes the network problem.

If $p_e(x(e)) >0$ for all $e$, then the network problem has unique solution. Therefore, equality
holds in (\ref{eqn:strictinequality}) if and only if $v=x$, i.e, $x$ is a fixed point of the
mapping
$T$. Thus we have
\begin{align}
0 &\leq \sum_{e=1}^{n} p_e(x(e)) \log~\frac{v(e)}{x(e)}\stackrel{(a)}{\leq} \sum_{e=1}^{n} p_e(x(e)) \left(\frac{v(e)}{x(e)}-1\right)\nonumber\\
&=\sum_{e=1}^{n}w_e^{\prime}(x(e))(v(e)-x(e))=\nabla W(x)\bigcdot (v-x), \label{eqn:timederivativepositive}
\end{align}
where (a) uses the inequality $\log\,y\leq y-1$.

More generally, let $p_e(x(e))>0$ for $e \in S \subset [n]$ and $p_e(x(e))=0$ for $e \in S^{c}$; in
particular, $x(e)=0$ for $e \in S^c$. Define $\tilde{v}$ to be $\tilde{v}(e)=v(e)1_{S}(e), e \in
[n]$.

The value of the objective function in (\ref{eqn:functionT}) evaluated at $v$ and $\tilde{v}$ are equal. Hence $\tilde{v} \in T(x)$,
\begin{align}\label{eqn:strictinequality1}
\sum_{e\in S}p_{e}(x(e))\log(\tilde{v}(e))\geq \sum_{e\in S}p_{e}(x(e))\log(x(e)),
\end{align}
and $\tilde{v}$ must be the unique solution to the problem defined in (\ref{eqn:subnetprob2}).
Therefore,
(\ref{eqn:strictinequality1}) holds with equality if and only if $\tilde{v}=x$.
Following the steps leading to (\ref{eqn:timederivativepositive}), we have
\begin{align}\label{eqn:strictinequality2}
\sum_{e\in S}w_e^{\prime}(x(e))(\tilde{v}(e)-x(e)) \geq 0
\end{align} which is a strict inequality if $\tilde{v}\neq x$.
Since $v(e)-x(e) \geq 0$ for $e \in S^{c}$, this along with (\ref{eqn:strictinequality2}) yields
\begin{align}\label{eqn:strictinequality3}
\nabla W(x)\bigcdot(v-x) \geq 0;
\end{align}
since $w_e^{\prime}(x(e))=w_e^{\prime}(0) > 0$ for $e \in S^c$, equality holds in (\ref{eqn:strictinequality3}) if and only if $v=x$. Hence
\begin{align}\label{eqn:rateofchange}
\frac{dV(x(t))}{dt}= -\nabla W(x(t))\bigcdot (v-x(t)) \leq 0 , ~ \forall ~ v \in T(x(t)).
\end{align}
The inequality in (\ref{eqn:rateofchange}) holds with an equality if and only if $v=x$. Therefore $V$ is a Lyapunov function for $\Lambda$. $\hfill \blacksquare$
\end{proof}

\begin{lemma}\label{lem:ictsubsetlambda}
Let $\Lambda$ be the set of fixed points of $T$. Every internally chain transitive set for $F$ in
(\ref{eqn:Fdefinition}) is a singleton that is a subset of $\Lambda$.
\end{lemma}
\begin{proof}
By Lemma \ref{lem:finitefixedpoint}, there are at most finitely many fixed points of the mapping
$T$.
Hence the  cardinality of the set $\Lambda$ is finite and $V(\Lambda)$ has empty interior. Also, by
Lemma \ref{lem:Vlyapunov}, $V$ is a Lyapunov function for $\Lambda$. Proposition 3.27 of
\cite{2005xxSIAM_BenHofSor}  states that if $V$ is a Lyapunov function for $\Lambda$ and if
$V(\Lambda )$ has an empty interior, then every internally chain transitive set is a subset of
$\Lambda$.

Choose $\epsilon$ small enough so that open balls of radius $\epsilon$ centered at each of the finite
points $\Lambda$ are disjoint. Fix $T \geq 0$. Since any $(\epsilon,T)$ chain involves remaining in
$\Lambda$ for all time and jumps of size at most $\epsilon$ to another point in $\Lambda$, by the
disjointedness of the $\epsilon$-balls covering $\Lambda$, there can be no $(\epsilon,T)$ chain
in $\Lambda$ joining two of its distinct points. It follows that the internally chain transitive
subsets of $\Lambda$ are singletons. $\hfill \blacksquare$
\end{proof}
\begin{lemma}
The iterates $x^{(k)}$ converges to a fixed point of the mapping $T$.
\end{lemma}
\begin{proof}
In Lemma \ref{lem:rxisps}, we showed that $r_x(t)$ is a perturbed solution to (\ref{eqn:DI}). By
Lemma \ref{lem:psisict}, the limit set of $r_x(t)$ is internally chain transitive. By Lemma
\ref{lem:ictsubsetlambda}, $L(r_x)$ is a singleton and $L(r_x) \subset \Lambda$. Let $\hat{x} \in
L(r_x)$. Since $\mathcal{A}$ is compact and $\hat{x}$ is the only limit point of the sequence
$x^{(k)}$, every
subsequence of $x^{(k)}$ has a further subsequence that converges to $\hat{x}$. Hence $x^{(k)}$
converges to $\hat{x}$. $\hfill \blacksquare$
\end{proof}
In the rest of this section, we show that the iterates converge to $x^{\star}$, the optimal
solution to the system problem.

Let the dual variables of the optimization problem Network$(A,c;p)$ be $\mu_l, l \in [m]$.
Kelly et al. \cite{1998xxJORS_KelMauTan} simplified the dual to this problem to be:
\begin{align}
&\text{Dual}(p,A,c): \label{eqn:dualT}\nonumber \\
&~~\min_{\mu_l \geq 0, l \in [m]} \left(\sum_{e=1}^{n}p_e\cdot \log\frac{1}{\sum_{l:l\in e}\mu_l} +
\sum_{l=1}^{m}\mu_lc(l) \right).
\end{align}
We now argue that the search for the optimal $\mu_l,l\in [m]$ may be restricted to a compact set.
\begin{lemma}\label{lem:dualcompact}
The optimization problem in (\ref{eqn:dualT}) with $p_e=p_e(x(e))$ is equivalent to the following
optimization problem. For any $x \in \mathcal{A}$,
\begin{align}\label{eqn:dualobjfun}
\max_{0\leq \mu_l \leq 2P/c(l)} \sum_{e=1}^{n}p_e(x(e))\cdot\log \left(\sum_{l:l\in e}\mu_l\right)-\sum_{l=1}^{m}\mu_lc(l),
\end{align} where $P\coloneqq \max_{x \in \mathcal{A}}\sum_{e=1}^{n} x(e)\cdot w_e^{\prime}(x(e)) < \infty.$
\end{lemma}
\begin{proof}
Define $R(\mu)$ to be the objective function in (\ref{eqn:dualobjfun}). For any $\mu_l > 2P/c(l)$, by reducing $\mu_{l}$, we increase the
objective function's value. To see this, it suffices to show that $\frac{\partial R(\mu)}{\partial \mu_l} < 0$ for any $\mu_l > 2P/c(l)$.
But this is easily checked as follows:
\begin{align}
\frac{\partial R(\mu)}{\partial \mu_l}&=\sum_{e:e\ni l}p_e(x(e))\frac{1}{\sum_{l^{\prime}:l^{\prime}\in e}\mu_{l^{\prime}}}-c(l)\leq \sum_{e:e \ni l} p_e(x(e)) \frac{1}{\mu_l}-c(l) \nonumber \\ &\leq \frac{1}{\mu_l} \sum_{e=1}^{n}p_e(x(e)) -c(l) \nonumber \leq
\frac{1}{\mu_l}\left[\max_{x \in \mathcal{A}}\sum_{e=1}^{n} p_e(x(e))\right]-c(l) \nonumber \\&= \frac{P}{\mu_l}-c(l) < 0,
\end{align}
where the last inequality follows if $\mu_l > 2P/c(l)$. $\hfill \blacksquare$
\end{proof}
\begin{lemma}\label{lem:convxstar}
Let $x^{(k)}$ converge to $\hat{x}$, a fixed point of the mapping $T$. Then $\hat{x}= x^{\star}$,
the optimal solution to the system problem.
\end{lemma}
\begin{proof}
Let $v^{(k+1)}$ solve problem (\ref{eqn:functionT}) with $p_e=p_e(x^{(k)}(e))$, and so
$v^{(k+1)} \in T(x^{(k)})$ satisfies the KKT conditions
\begin{align}
&\frac{p_e(x^{(k)}(e))}{v^{(k+1)}(e)} - \sum_{l:l \in e} \mu_{l}^{(k)} +\eta_e^{(k)}=0, ~ e \in [n],\label{eqn:lagderiv} \\
&\mu_l^{(k)}\cdot \left( \sum_{e:e \ni l}v^{(k+1)}(e)-c(l)\right) = 0,~ l \in [m], \label{eqn:compslackintermediate}\\
&\eta_e^{(k)}\cdot v^{(k+1)}(e) = 0, ~e \in [n], \label{eqn:compslackintermediate2}\\
&\eta_e^{(k)}\geq 0,~ e \in [n],~ \mu_l^{(k)} \geq 0 , ~ l \in
[m].\label{eqn:dualpositiveintermediate}
\end{align}
Let us first claim that $x^{(k)} > 0$ for all $k \geq 0$. This holds for $k=0$, the initial point,
in Algorithm \ref{alg:mainalgorithm}. If, for some $k$, $x^{(k)}>0$, then $p_e(x^{(k)}(e)) > 0$ for
all $e$ and so, $v^{(k+1)} > 0$ and consequently, $x^{(k+1)}$ being a convex combination of
$x^{(k)}$ and $v^{(k+1)}$ also satisfies $x^{(k+1)} >0$. The claim follows by induction.
Since $x^{(k)}(e) > 0$, we have $v^{(k+1)}(e) > 0$, and so, by (\ref{eqn:compslackintermediate2}),
$\eta_e^{(k)}=0, ~ e \in [n]$.
Thus (\ref{eqn:lagderiv}) simplifies to
\begin{align}\label{eqn:denmupositive}
\frac{p_e(x^{(k)}(e))}{v^{(k+1)}(e))}=\sum_{l:l \in e} \mu_l^{(k)}.
\end{align}
Since we also have $p_e(x^{(k)}(e)) > 0$, and $v^{(k+1)}(e) > 0$ in (\ref{eqn:denmupositive}), we
have $\sum_{l:l \in e} \mu_l^{(k)} >0$. Hence
\begin{align}
v^{(k+1)}(e) &= \frac{p_e(x^{(k)}(e))}{\sum_{l:l \in e} \mu_l^{(k)}}=\frac{x^{(k)}(e)\cdot w_e^{\prime}(x^{(k)}(e))}{\sum_{l:l \in e} \mu_l^{(k)}}.\label{eqn:tautoinfty}
\end{align}

Suppose $\hat{x}(e)=0$. Observe that $x^{(k)}(e) > 0$ for all $k$. Hence $v^{(k+1)}(e) < x^{(k)}(e)$ infinitely often, which is the same as
saying
\begin{align}
\frac{x^{(k)}(e)w_e^{\prime}(x^{(k)}(e))}{\sum_{l:l \in e} \mu_l^{(k)}} < x^{(k)}(e)
\end{align}
occurs infinitely often.
This implies
\begin{align}\label{eqn:strictdecrease}
w_e^{\prime}(x^{(k)}(e)) < \sum_{l:l \in e} \mu_l^{(k)}
\end{align}
infinitely often.
Consider the subsequence that satisfies (\ref{eqn:strictdecrease}). Henceforth, let $x^{(k)}$ denote
that subsequence. We now make the following observations. In Lemma \ref{lem:dualcompact}, we showed
 that $\mu_l^{(k)}$ takes values on a compact set, and so we can find a further subsequence such
that $\mu_l^{(k_\tau)}\rightarrow \mu_l^{\star}$ for some $\mu_l^{\star}$, but for all $l \in [m]$.
Since (\ref{eqn:strictdecrease}) holds for $k = k_\tau$, by letting $\tau \rightarrow \infty$, we have
\begin{align}
w_e^{\prime}(\hat{x}(e)) &\leq \sum_{l:l \in e} \mu_l^{\star}, ~ \forall ~ e~\text{such that}~ \hat{x}(e)=0 \label{eqn:lagderivativefinal}\\
w_e^{\prime}(\hat{x}(e)) &= \sum_{l:l \in e} \mu_l^{\star},~ \forall ~ e~\text{such that} ~\hat{x}(e)>0;
\end{align}
the latter follows by letting $\tau \rightarrow \infty$ in (\ref{eqn:tautoinfty}) and from
$x^{(k+1)}=x^{(k)}+a_k (v^{(k+1)}-x^{(k)})$.
Choose
\begin{align}
\eta^{\star}(e)& =0 ~\text{if} ~\hat{x}(e) > 0, \\
\eta^{\star}(e)&=\sum_{l:l \in e} \mu_l^{\star}-w_e^{\prime}(\hat{x}(e))~\text{if} ~\hat{x}(e) = 0.
\end{align}
Thus, from (\ref{eqn:lagderivativefinal}),
\begin{align}
\eta^{\star}(e)\cdot \hat{x}(e)=0~\forall ~ e \in [n] ~\text{and}~\eta^{\star}(e) \geq 0~ \forall ~
e \in [n].
\end{align}
Since (\ref{eqn:compslackintermediate}) and (\ref{eqn:dualpositiveintermediate}) are true for
indices $k_\tau$, taking limit as $\tau \rightarrow \infty$, we get
\begin{align}
\mu_l^{\star}\cdot \left( \sum_{e:e \ni l}\hat{x}(e)-c(l)\right)& = 0,~ l\in [m],\\
\mu_l^{\star} &\geq 0, ~ \forall ~l.\label{eqn:positivityfinal}
\end{align}
Equations (\ref{eqn:lagderivativefinal})-(\ref{eqn:positivityfinal}) are the KKT conditions for the system problem. Hence $\hat{x}=x^{\star}$,
the optimal solution to the system problem. $\hfill \blacksquare$
\end{proof}

\vspace*{.25cm}

\noindent {\em Proof of Theorem \ref{thm:maintheorem}:}
Theorem \ref{thm:maintheorem} follows from Lemma \ref{lem:convxstar}. $\hfill \blacksquare$

\section{Algorithmic Complexity and Speed of Convergence}
\label{sec:complexity}
In this section, we remark on the complexity of Algorithm \ref{alg:mainalgorithm}. Each iteration
of
the algorithm has 1) a user update which adapts the amount a user is willing to pay to the
network,  and 2) a network update which adapts the rates allocated to the users.

Since $w_e$ is known at the user end, $w_e^{\prime}$ is easy to obtain either numerically or
analytically. Hence the user update (\ref{eqn:userupdate}) can be implemented by each user in
$\mathcal{O}(1)$ steps.

\begin{figure}[htb]
\centering
\includegraphics[scale=0.6]{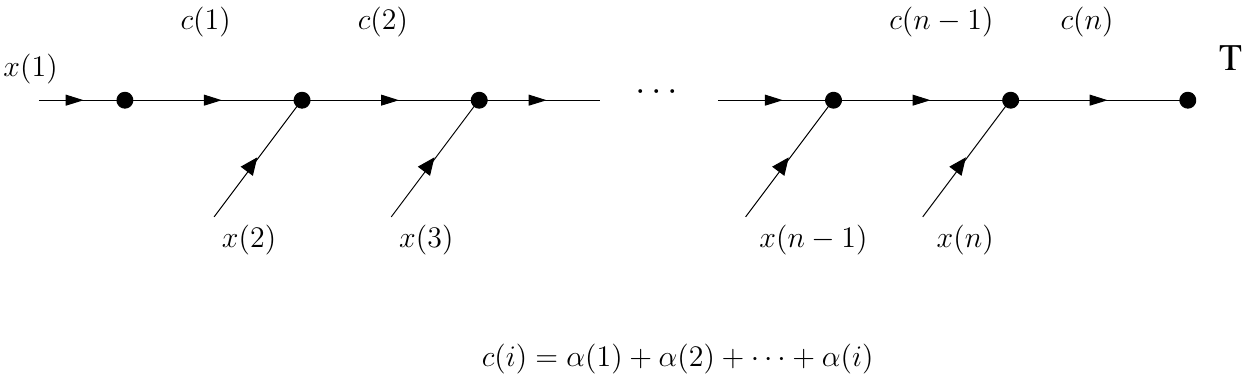}
 \caption{Flow aggregating network.}
 \label{fig:figure1}
 \end{figure}

The network update consists of solving the network problem (\ref{eqn:functionT}). Its complexity
depends on the network structure. We indicate the complexity of the network update for the
following simple networks: a polymatroidal network with a single source and multiple sinks or
multiple sources and a single sink; a flow aggregating network with the structure in Figure
\ref{fig:figure1}.

{\em Polymatroidal network}: Consider a network with a single source and $n$ sinks. The
source sends flows at a rate $x(1),x(2),\ldots,x(n)$ to the sinks. Megiddo \cite{1974xx_Meg} showed
that the set of feasible flows $(x(e), e \in [n])$ forms the independent set of a polymatroid.
Therefore, the network problem is a separable concave maximization over the independent set of a
polymatroid. The fastest known algorithm that solves this optimization problem is a scaling based
greedy algorithm proposed by Hochbaum \cite{1994xxMOR_Hoc}. The algorithm obtains an
``$\epsilon$-optimal" solution to the network problem in $\mathcal{O}\left(n (\log n +F) \log
\frac{B}{n\epsilon}\right)$ where $F$ is the complexity to check whether a certain increase in one
of the  components of $(x(e), e \in [n])$ would make the flow infeasible. $B$ is the total amount of
resource to be allocated and is $\mathcal{O}(n)$.

{\em Flow aggregating network}:
Let  $(x(e),e \in [n])$ denote the flow through the network in Figure \ref{fig:figure1}. The flow
constraints of the network are
\begin{align}
\begin{split}\label{eqn:lac}
x(e) & \geq 0, e \in [n],\\
x(1) &\leq \alpha(1),\\
x(1)+x(2) &\leq \alpha(1)+\alpha(2),\\
&~\,\vdots\\
x(1)+x(2)+\cdots+x(n) &\leq \alpha(1)+\alpha(2)+\cdots+\alpha(n).
\end{split}
\end{align}
where $\alpha(e) \geq 0$ for all $e$.
The constraints in (\ref{eqn:lac}) are referred to as linear ascending constraints. This problem
arises as the core optimization problem in several wireless communication problems (Padakandla and
Sundaresan \cite{200910TCOM_PadSun}, Viswanath and Anantharam \cite{200206TIT_VisAna}, Lagunas
et al. \cite{200405TSP_PalLagCio}, Sanguinetti et al. \cite{201205TSP_SanDam}) and operations
research problems (Clark and Scarf \cite{1960xxMS_ClaSca}, Wang \cite{201510OL_Ziz}). See
\cite{2016xx_AkhSun} for a survey and a discussion of several algorithms. The network problem is
the maximization of a so-called $d$-separable concave function over the linear ascending
constraints. Veinott Jr. \cite{1971xxMS_Vei} mapped this problem to the geometrical problem of
finding the ``concave cover" of the set of points $(\sum_{e=1}^{i}\alpha(e),\sum_{e=1}^{i}p_e),~ i
\in [n]$, in $\mathbb{R}^2$. The ``string algorithm" of Muckstadt and Sapra \cite{2010xxSPR_MucSap}
finds the concave cover of a set of points in $\mathbb{R}^{2}$ in $\mathcal{O}(n)$ steps. See
\cite{2016xx_AkhSun} for details.

{\em Simulations}:
We now discuss some simulation studies investigating the speed of convergence of the ODE that our proposed algorithm will track. But we caution the reader that the ODE convergence rate does not give the full picture of convergence rate since the timescale is dictated by the step sizes. In the plots, the solid curves correspond to the error plots of the system
\begin{align}
 \label{eqn:scaledmainalgdi} \frac{d}{dt}x(t)  \in \kappa \cdot \left ( T(x(t)) -x(t) \right ).
\end{align}
The dashed curves correspond to the error plots of the system (\ref{eqn:price})-(\ref{eqn:cost}).
The differential inclusion (\ref{eqn:scaledmainalgdi}) has scaling factor $\kappa$ when compared
with (\ref{eqn:mainalgdi}) and corresponds to a scaled version of Algorithm \ref{alg:mainalgorithm}.
The scaling is to enable comparison of (\ref{eqn:scaledmainalgdi}) with the system
(\ref{eqn:price})-(\ref{eqn:cost}) which already has the scaling factor  $\kappa$ in
(\ref{eqn:update}).

\begin{figure}[t]
\centering
\includegraphics[scale=0.5]{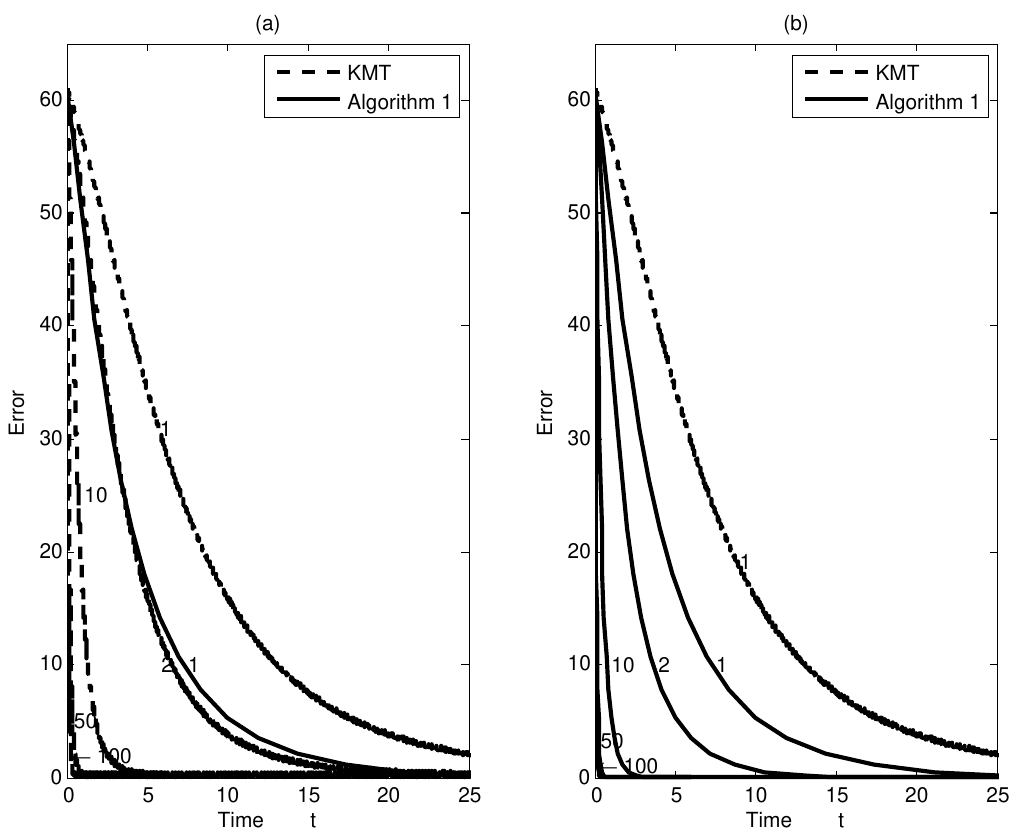}
 \caption{(a) Error dynamics for (\ref{eqn:price})-(\ref{eqn:cost}) denoted (``KMT'') with
$\kappa = 1,2,10,50,100$; error dynamics for (\ref{eqn:scaledmainalgdi}) denoted ``Algorithm
\ref{alg:mainalgorithm}'' with $\kappa=1$ is also plotted for comparison. (b) The roles of
(\ref{eqn:price})-(\ref{eqn:cost}) and (\ref{eqn:scaledmainalgdi}) are swapped.}
 \label{fig:figure_kmt_alg1}
 \end{figure}

 \begin{figure}[hbt]
\centering
\includegraphics[scale=0.7]{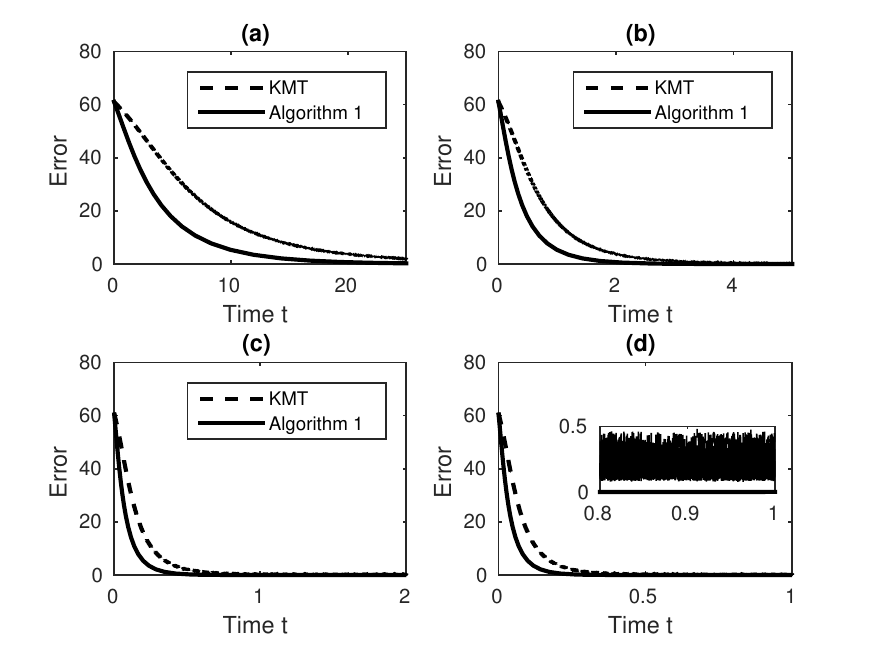}
 \caption{Comparison of error dynamics for the systems (\ref{eqn:price})-(\ref{eqn:cost}) (``KMT''
curve) and (\ref{eqn:scaledmainalgdi}) (``Algorithm \ref{alg:mainalgorithm}'' curve). Subplot (a)
has $\kappa=1$, (b) has $\kappa=10$, (c) has $\kappa=50$, and (d) has $\kappa=100$. Note that the
time axis is scaled differently in each of the subfigures. The inset picture in (d) shows an enlarged view of the error. The ``KMT'' curve settles at a small but positive error.}
 \label{fig:figure_kmt_alg1_identical_kappa}
 \end{figure}

All figures are for the flow aggregating network with $n=m=10$ and link capacities $c(l)=10 \times l, l \in [m]$. The utility functions are chosen as $(1/\beta (e))\cdot x^{\beta (e)}$, for $ e \in
[n]$ with $\beta(e)=.09\cdot e$. The initial point for both differential inclusions is always the lexicographically maximal point\footnote{The lexicographically maximal point is one where the minimum allocation (across users) is maximized among all feasible points; further the second minimum is maximized among all points with equal minimum allocation, and so on.}. This is the most natural starting point when the network does not know the users' utility functions and considers all users to be equal. While we report the results only for this particular $\beta$ and $c$, we have simulated several other settings, and the results are qualitatively the same. We do not repeat them here for brevity.

Figure \ref{fig:figure_kmt_alg1}(a) shows that as $\kappa$ scales up, the speed of convergence of the
system (\ref{eqn:price})-(\ref{eqn:cost}) increases. For comparison, we have included the solid
curve for (\ref{eqn:scaledmainalgdi}) with $\kappa=1$.

Figure \ref{fig:figure_kmt_alg1}(b) shows that as $\kappa$ scales up, the rate of convergence of
(\ref{eqn:scaledmainalgdi}) also increases similarly. Again, for comparison, we have included the
dashed curve for (\ref{eqn:price})-(\ref{eqn:cost}) with $\kappa=1$.

These two subfigures show that convergence can be sped up similarly in the two systems,
(\ref{eqn:price})-(\ref{eqn:cost}) and (\ref{eqn:scaledmainalgdi}), by simply increasing $\kappa$.

Figures \ref{fig:figure_kmt_alg1_identical_kappa}(a)-\ref{fig:figure_kmt_alg1_identical_kappa}(d) compare (\ref{eqn:price})-(\ref{eqn:cost})
directly with (\ref{eqn:scaledmainalgdi}) for identical $\kappa$. The speeding up parameter
$\kappa$ equals 1, 10, 50 and 100 in Figures
\ref{fig:figure_kmt_alg1_identical_kappa}(a), \ref{fig:figure_kmt_alg1_identical_kappa}(b), \ref{fig:figure_kmt_alg1_identical_kappa}(c), and \ref{fig:figure_kmt_alg1_identical_kappa}(d)
respectively. These figures demonstrate that convergence speed of (\ref{eqn:scaledmainalgdi}) is comparable to that of
(\ref{eqn:price})-(\ref{eqn:cost}) as long $\kappa$ is identical for the two
systems. We saw the same qualitative behavior across several randomly chosen problem parameters.

The inset in Figure \ref{fig:figure_kmt_alg1_identical_kappa}(d) shows an enlarged view of the error plots after the algorithms' settlement close to their respective limiting values. We see that the error plot of system (\ref{eqn:price})-(\ref{eqn:cost}) settles at a small but positive value. This is consistent with the observation that KMT algorithm solves a relaxation of the original system problem.

\begin{figure}[t]
\centering
\includegraphics[scale=0.8]{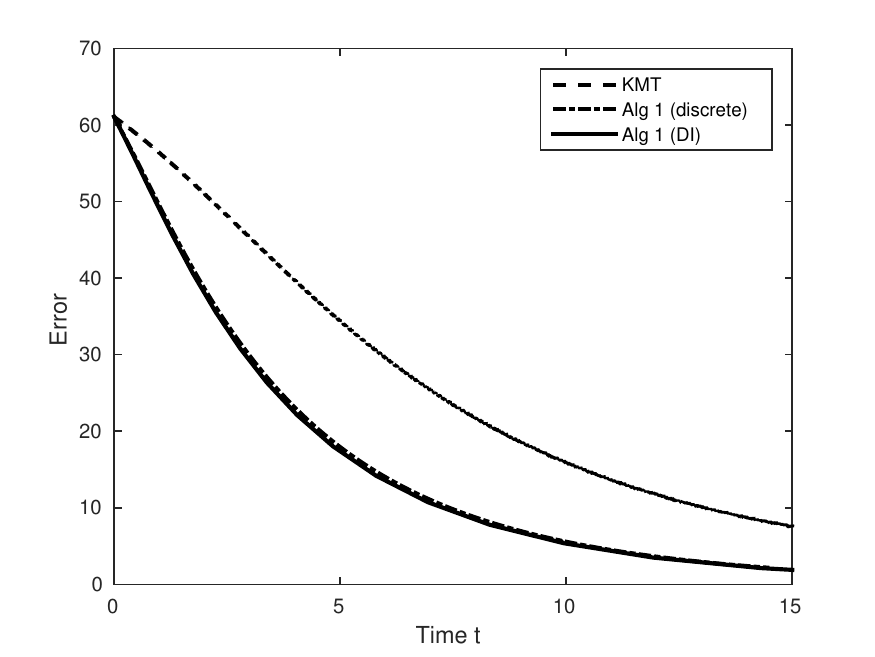}
 \caption{This figure is generated using the same problem parameters used to generate Figure \ref{fig:figure_kmt_alg1_identical_kappa}(a). The curve KMT corresponds to system (\ref{eqn:price})-(\ref{eqn:cost}) and curve Alg1 (DI) plots the ODE (\ref{eqn:scaledmainalgdi}) for $\kappa=1$. Alg1 (discrete) plots the iterates of Algorithm \ref{alg:mainalgorithm} at timeinstants $t_k$ given by (\ref{eqn:comp_dis_cts}).}
 \label{fig:figure_discrete_alg1}
 \end{figure}

Figure \ref{fig:figure_discrete_alg1} reproduces the plots in Figure \ref{fig:figure_kmt_alg1_identical_kappa}(a) but with abscissa values restricted to time interval $[0,15]$. The dash-dotted line plots the iterates put out by Algorithm \ref{alg:mainalgorithm} with the $k^{th}$ iterate plotted at time instant
\begin{align}\label{eqn:comp_dis_cts}
t_k=\sum_{i=1}^{k}a_i.
\end{align}
As expected, we see that the iterates trace the error plot of the ODE (\ref{eqn:scaledmainalgdi}) for $\kappa=1$. This justifies the comparison between the system (\ref{eqn:price})-(\ref{eqn:cost}) and the ODE (\ref{eqn:scaledmainalgdi}).

 \begin{figure}[t]
\centering
\includegraphics[width=12cm, height =7cm]{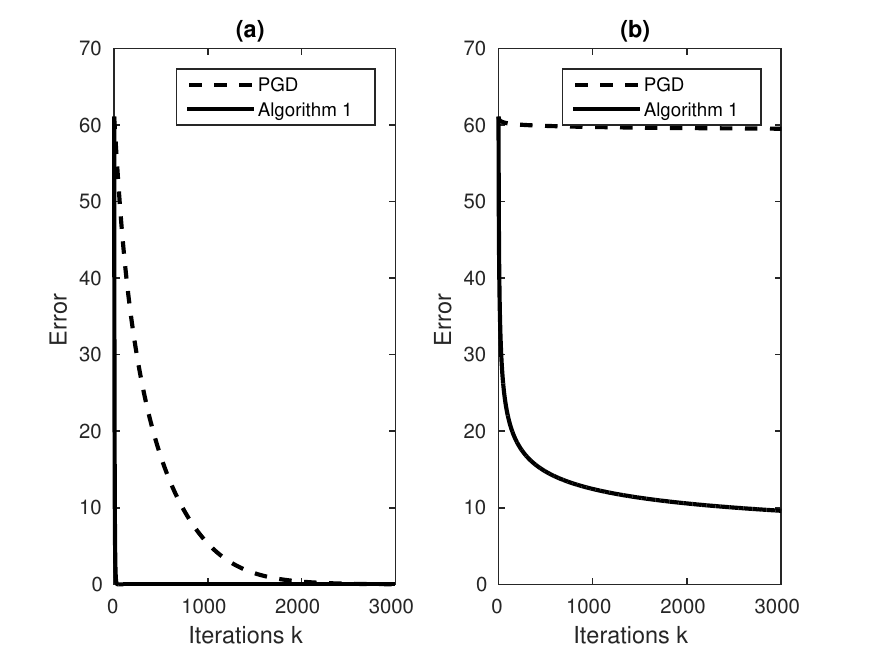}
 \caption{Figure compares the performance of Algorithm \ref{alg:mainalgorithm} with the projected gradient descent algorithm. The plot in Figure \ref{fig:figure_pgd_alg1}(a) corresponds to the case when step-size for both algorithms is calculated using Armijo-Goldstein rule whereas for the plots in Figure \ref{fig:figure_pgd_alg1}(b) step-size $a_k=\frac{1}{k+1}$ is used.}
 \label{fig:figure_pgd_alg1}
 \end{figure}

We also compare our algorithm with a benchmark interior point algorithm, the projected gradient algorithm \cite[Sec. 2.3]{bertsekas1999nonlinear}. The projected gradient descent algorithm is not distributed because the step-size selection according to Armijo-Goldstein rule would require the knowledge of utility functions. Hence the comparison is made in the following two ways.

We first compare the case when the stepsize is $a_k = 1/(k+1)$ as for stochastic approximation. With these fixed stepsizes, the projected gradient descent can also be implemented in a distributed fashion, similar to ours. The network asks all users to send flows according to $x^{(k)}$ and invites these users to send gradients of their private utility functions at these points. The users follow this. With the gradient information, the network identifies a new location by employing gradient descent, projects it on the feasible set, and then asks users to send flows according to this projected $x^{(k+1)}$. The procedure then repeats.

We next compare the case when the stepsizes are according to the Armijo-Goldstein rule. This cannot be done in a distributed fashion since the improvement comparisons require knowledge of the private utility functions. So, for fair comparison, we too use stepsizes according to the Armijo-Goldstein rule to get a centralized variant of Algorithm 1.

As can be seen from the two new plots in Figures \ref{fig:figure_pgd_alg1}(a) and \ref{fig:figure_pgd_alg1}(b), in both cases, our algorithm does much better than projected gradient descent. In Figure \ref{fig:figure_pgd_alg1}(b), our algorithm is very close to the axes. This is quite reassuring.

\section{Conclusion}
\label{sec:conclusion}
We considered the network utility maximization problem in a distributed
framework where the users do not know the network structure or utility functions of other users and
the network does not know the users' utility functions. We decomposed the system problem into user
subproblems and a network subproblem following the methodology of \cite{1998xxJORS_KelMauTan}.
Unlike the dual decomposition iterative methods of \cite{1958xxEJES_ArrHur},
\cite{1998xxJORS_KelMauTan}, \cite{1999xx_LowLap}, etc., the iterations proposed in Algorithm
\ref{alg:mainalgorithm} ensure feasibility at every step. The convergence of the algorithm was shown using the theory of differential inclusions. The iterates avoid local maxima traps on the facets. Efficient methods to solve the network problem for some
special networks were also described. Finally, sample simulations show that, in several examples,
Algorithm \ref{alg:mainalgorithm}'s associated differential inclusion (\ref{eqn:mainalgdi})
converges faster to the system optimal point when compared with the iterates arising from the ODE
(\ref{eqn:price})-(\ref{eqn:cost}). The ODE convergence rate however does not give the full picture of convergence rate of the iterates since the timescale is dictated by the step sizes. For the convergence rate of the iterates, a natural approach is to use the method of Borkar \cite[Ch. 4]{2008xxCUP_Bor} to get sample complexity bounds. However they do not directly apply since the ODE dynamics is not necessarily Lipschitz, which is a crucial assumption in \cite[Ch.~4]{2008xxCUP_Bor}. See \ref{appsec:DiscontinuousT} for a discontinuous $T$ mapping. A more intricate analysis of convergence rates is therefore required and is left as future research.

\appendix
\section{Counterexample to the Algorithm of Hou et al. \cite{2010xxINF_HonKum}}\label{sec:stepsizecounterexample}
In this section, we provide an example where Algorithm \ref{alg:mainalgorithm} does not converge for a constant step size $a_k=\gamma,~ \gamma \in (0,1) $. Consider a two user single link network. The system problem
for the case is
\begin{align}
 \text{Maximize}
 & \hspace{3mm}w_1(x(1))+w_2(x(2)) \nonumber\\
 \text{subject to}&\hspace{3mm}x(1)+x(2)=B, \nonumber
 \end{align}
 where $B$ is the capacity of the link. Let $d=(d_1,d_2)$ be the initial flow through the link. Let $f=(f_1,f_2)$ be defined as
 \begin{align}
  f=d+\gamma (T(d)-d),\label{eqn:firstupdate}
 \end{align}
the flow allocated in the first iteration of Algorithm \ref{alg:mainalgorithm}. Without loss of generality, choose
\begin{align}\label{eqn:flow1condition}
 d_1 < f_1 < T_1(d) < B.
\end{align}
Also, choose
\begin{align}\label{eqn:equaljumps}
 T_1(f)=d_1+f_1-T_1(d)
\end{align}
The flow allocated to user 1 in the second iteration is
\begin{align}
 f_1+\gamma (T_1(f)-f_1)&\stackrel{(a)}{=}f_1 + \gamma (d_1-T_1(d)) \nonumber \\
                       &\stackrel{(b)}{=} d_1 +\gamma (T_1(d)-d) + \gamma (d-T_1(d))\nonumber\\
                       &=d_1,\label{eqn:oscillation}
\end{align}
where (a) and (b) are due to (\ref{eqn:equaljumps}) and (\ref{eqn:firstupdate}) respectively. Equations (\ref{eqn:firstupdate}) and (\ref{eqn:oscillation}) imply that the
flows put out by the algorithm oscillates from $d$ to $f$ and vice versa. It remains to be shown that there exists $w_1$ and $w_2$ that is
consistent with the choices made in (\ref{eqn:flow1condition}) and (\ref{eqn:equaljumps}).

We have, by the definition of $T(.)$,
\begin{align}\label{eqn:T1d}
 \frac{d_1w_1^{\prime}(d_1)}{d_1w_1^{\prime}(d_1)+d_2w_2^{\prime}(d_2)}=\frac{T_1(d)}{B},\\
 \frac{f_1w_1^{\prime}(f_1)}{f_1w_1^{\prime}(f_1)+f_2w_2^{\prime}(f_2)}=\frac{T_1(f)}{B}.\label{eqn:T1f}
\end{align}
We rewrite the equations (\ref{eqn:T1d}) and (\ref{eqn:T1f}) as
\begin{align}
 d_1(1-\frac{T_1(d)}{B})w_1^{\prime}(d_1)-d_2\frac{T_1(d)}{B}w_2^{\prime}(d_2)=0, \label{eqn:l1}\\
  f_1(1-\frac{T_1(f)}{B})w_1^{\prime}(f_1)-f_2\frac{T_1(f)}{B}w_2^{\prime}(f_2)=0. \label{eqn:l2}
\end{align}

\begin{figure}
\centering
\includegraphics[draft=false,scale=0.35]{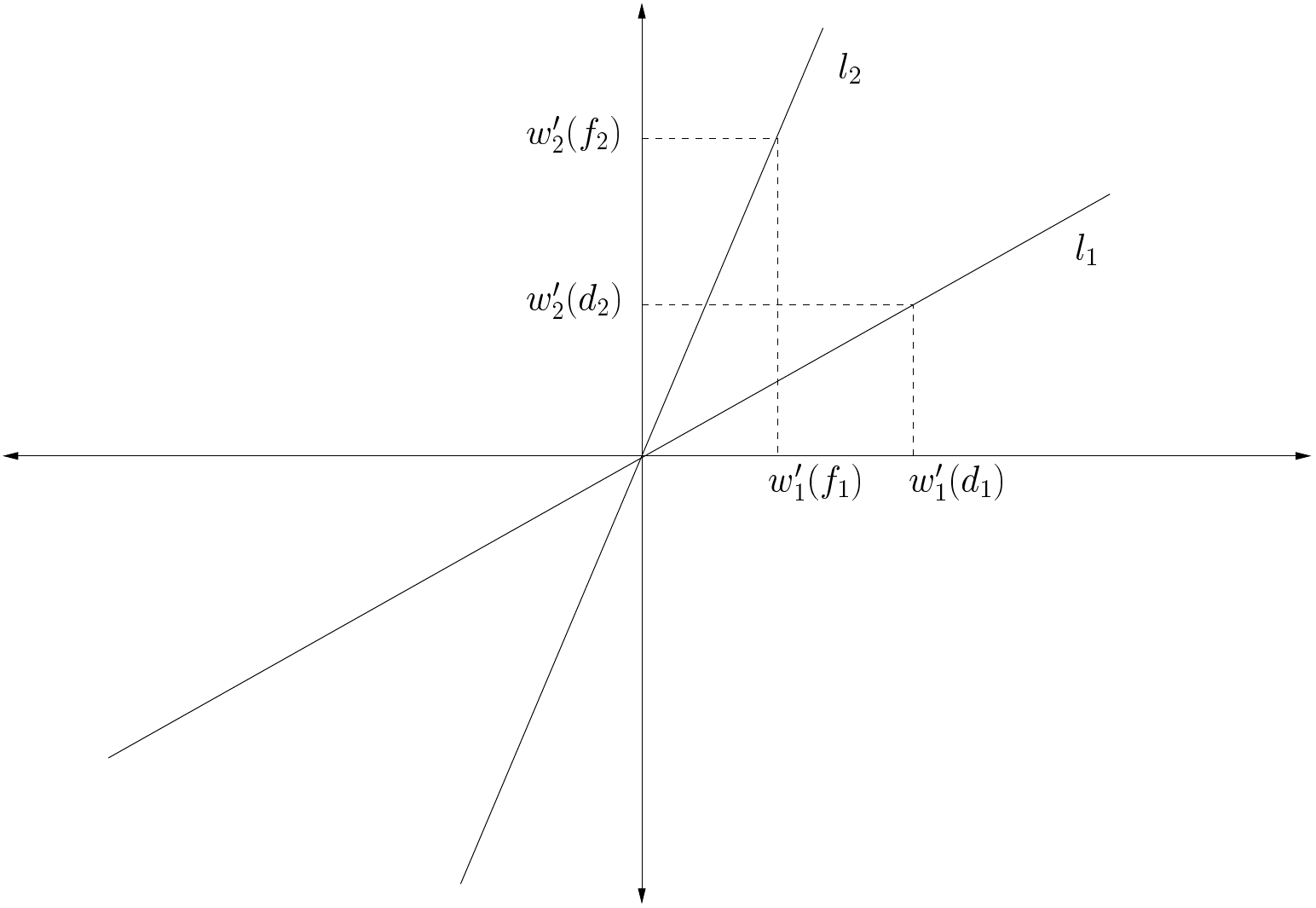}
\caption{Figure shows lines $l_1$ and $l_2$. $w_1^{\prime}(d_1)$ and $w_1^{\prime}(f_1)$ are plotted on the horizontal axis.
$w_2^{\prime}(d_2)$ and  $w_2^{\prime}(f_2)$ are plotted on the vertical axis.}
\label{fig:figurecounter}
\end{figure}

We now view (\ref{eqn:l1}) and (\ref{eqn:l2}) as linear equations in $w_1^{\prime}(d_1),w_2^{\prime}(d_2)$ and
$w_1^{\prime}(f_1),w_2^{\prime}(f_2)$ respectively. Lines $l_1$ and $l_2$ in Figure \ref{fig:figurecounter} plot (\ref{eqn:l1}) and
(\ref{eqn:l2}) respectively. Since $d_1+d_2=f_1+f_2=B$, $d_1 <f_1$ implies $d_2 >f_2$. Also, from (\ref{eqn:oscillation}) and the fact that
$d_1<f_1$, we have $T_1(f) < d_1$. Hence, by (\ref{eqn:flow1condition}), $T_1(f) < T_1(d)$.  Since $d_1 <f_1$, $d_2 > f_2$ and $T_1(f) <
T_1(d)$,
the slope of $l_1$ is smaller than the slope of $l_2$.

Since $d_1 <f_1$ and $d_2 > f_2$, by the strict concavity of $w_1$ and $w_2$, we must have
\begin{align}
 w_1^{\prime}(f_1) < w_1^{\prime}(d_1) \text{~and~} w_2^{\prime}(d_2) < w_2^{\prime}(f_2).
\end{align}

Figure \ref{fig:figurecounter} shows how to choose $w_1^{\prime}(d_1),w_2^{\prime}(d_2),w_1^{\prime}(f_1)$ and $w_2^{\prime}(f_2)$ satisfying
(\ref{eqn:flow1condition}),(\ref{eqn:equaljumps}),(\ref{eqn:l1}) and (\ref{eqn:l2}).

\section{An Example of a Discontinuous $T$ mapping}
\label{appsec:DiscontinuousT}
In this Appendix, we show that there is no selection from within $T(x)$ that could make the
selection a single continuous mapping. Consider a special case of $T(x)$ as defined below. Take
$w_e(\cdot)=w(\cdot)$ for some increasing and strictly concave $w(\cdot)$. Let
\begin{align}
T(x)=  \arg \max_y\sum_{e=1}^{3}p_e(x(e))\cdot &\log(y(e))\label{eqn:example}\\
\mbox{subject to    }~~~~~~~y(1) &\leq c,  \nonumber \\
                  y(1)+y(2) &\leq  2c,\nonumber\\
                   y(1)+y(2)+y(3) &\leq 3c. \nonumber
\end{align}
where $p_e(x(e))= x(e)\cdot w^{\prime}(x(e)),~e=1,2,3$. Let $w(\cdot)$ satisfy $x(e)\cdot
w^{\prime}(x(e))|_{x(e)=0}=0$. Consider $T(x)$ at $x= (c,0,0)$. We have
 \begin{align}
T(c,0,0)=\{(c,r,s):r\in [0,c], s \in [0,2c]:r+s \leq 2c\}. \nonumber
\end{align}
Consider  a sequence $y^{(k)} \rightarrow x$ such that $y^{(k)}(1)=c,~y^{(k)}(2)=0, ~ y^{(k)}(3)>0$
for each $k$ and $\lim_{k \rightarrow \infty}y^{(k)}(3)=0$. It is easy to see that
$T(y^{(k)})=\{(c,0,2c)\}$ for each $k$, and so we must select $(c,0,2c)$ at $x=(c,0,0)$.

Now, consider another sequence $z^{(k)}\rightarrow x$. Let $z^{(k)}(1)=c,~z^{(k)}(2)=z^{(k)}(3)> 0$
for each $k$ and $\lim_{k \rightarrow \infty}z^{(k)}(2)= \lim_{k \rightarrow \infty}z^{(k)}(3)=0$.
We then have $T(z^{(k)})=\{(c,c,c)\}$ for each $k$, and so we must now select $(c,c,c)$ at
$x=(c,0,0)$.
Since these two selections do not match, $T$ cannot be made continuous at $(c,0,0)$ by a choice of
a value in $T(c,0,0)$.

So $T$ has to be dealt with as a set valued mapping, which brings us to differential inclusions.

\bibliographystyle{IEEEtran}
\bibliography{IEEEabrv,wisl}

%







\end{document}